\documentclass[11pt]{article}
\usepackage{mathrsfs}
\usepackage{amsfonts}
\usepackage{indentfirst}
\usepackage{amstext, graphicx, amsmath, latexsym, amssymb, amscd, amsfonts}
\usepackage{color}
\usepackage{subfigure}
\usepackage{graphics}
\usepackage{graphicx}
\usepackage{epstopdf}
\usepackage{array}
\allowdisplaybreaks
\numberwithin{equation}{section}
\newtheorem{thm}{\bf Theorem}       
\newtheorem{cor}{\bf Corollary} 
\newtheorem{lem}{\bf Lemma}

\newtheorem{defn}{\bf Definition}

\numberwithin{lem}{section} \numberwithin{thm}{section}
\numberwithin{cor}{section}
\numberwithin{prop}{section}
\numberwithin{defn}{section}
\numberwithin{exa}{section} 
\begin{document}
	\title{\bf \Large Perturbation bounds for the largest $C$-eigenvalue of piezoelectric-type tensors}
	\author{Xifu Liu${}^{\mbox{\tiny a}}$,~Dongdong Liu${}^{\mbox{\tiny b}}$\thanks{%
			Corresponding author. \emph {E-mail addresses}: lxf211@cqnu.edu.cn (X. Liu), ddliu@gdut.edu.cn (D. Liu).},~Yaping Shi\,${}^{\mbox{\tiny a}}$  \\
\\
\small\textit{${}^{\mbox{\tiny a}}$School of Mathematical Sciences, Chongqing Normal University, Chongqing, China}\\
\small\textit{${}^{\mbox{\tiny b}}$School of Mathematics and Statistics, Guangdong University of Technology, Guangzhou, China}}

\date{}
\maketitle

\begin{abstract}
In this paper, we focus on the perturbation analysis of the largest $C$-eigenvalue of the piezoelectric-type tensor which has concrete physical meaning which determines the highest piezoelectric coupling constant. Three perturbation bounds are presented, theoretical analysis and numerical examples show that the third perturbation bound has high accuracy when the norm of the perturbation tensor is small.    
\end{abstract}

\section{Introduction}\label{s:1}

Let $\mathbb{R}$ be the real field. For a positive integer $n$, we denote the set $\{ 1,2, \cdots ,n\}$ by $[n]$. We call $\mathcal{A}=(a_{i_1 i_2 \cdots i_m})$ a real tensor of $m$-th order  $n$ dimension, denoted by $\mathcal{A} \in\mathbb{R}^{[m, n]}$ (or $\mathbb{R}^{n  \times n  \times \cdots \times n}$), if $a_{i_1 i_2 \cdots i_m} \in\mathbb{R}$ for all $i_1, \ldots, i_m \in [n]$. For any tensor $\mathcal{A} = (a_{i_1 i_2 \cdots i_m}) \in\mathbb{R}^{[m, n]}$, if its entries $a_{i_1 i_2 \cdots i_m}$ are invariant under any permutation of its indices, then $\mathcal{A}$ is called a symmetric tensor.

\par
Among various concepts of tensor eigenvalues, $Z$-eigenvalue [\ref{Qi1}] is one of the most popular research topics in tensor field.
Let $\mathcal{A} \in\mathbb{R}^{[m, n]}$. If there exist a complex number $\lambda$ and a nonzero complex vector $\mathbf{x} = (x_1, x_2, \cdots, x_n)^{T}$ that are solutions of the following polynomial equations,
$$\mathcal{A} \mathbf{x}^{m-1} = \lambda \mathbf{x},\;\;\mathbf{x}^{T} \mathbf{x}=1,$$
where $\mathcal{A} \mathbf{x}^{m-1}$ is a vector, and its $i$-th entry is given by
\begin{equation}\label{defAxm-1}
(\mathcal{A} \mathbf{x}^{m-1})_i = \sum\limits_{i_2, \cdots, i_m \in [n]} {a_{i i_2 \cdots i_m} x_{i_2} \cdots x_{i_m}}, i \in [n],
\end{equation}
then $\lambda$ is called an $E$-eigenvalue of $\mathcal{A}$ and $\mathbf{x}$ the $E$-eigenvector of $\mathcal{A}$ associated with $\lambda$, and $(\lambda, \mathbf{x})$ is called an $E$-eigenpair of $\mathcal{A}$.
We call $(\lambda, \mathbf{x})$ a $Z$-eigenpair if both $\lambda$ and $\mathbf{x}$ are real. Moreover, we have the following results on $Z$-eigenvalue.

\begin{lem}\label{l0}
$([\ref{Qi}])$ Let $\mathcal{A} \in\mathbb{R}^{[m, n]}$ be a symmetric tensor. Then $\mathcal{A}$ always has $Z$-eigenvalues. In particular, we have \\
$$\lambda_{Z \mathop {\max }}(\mathcal{A}) = \max \{ \mathcal{A} \mathbf{x}^m:\;  \mathbf{x}^T \mathbf{x} = 1,\; \mathbf{x} \in \mathbb{R}^n \},$$
$$\lambda_{Z \mathop {\min }}(\mathcal{A}) = \min \{ \mathcal{A} \mathbf{x}^m:\;  \mathbf{x}^T \mathbf{x} = 1,\; \mathbf{x} \in \mathbb{R}^n \},$$
where $\mathcal{A} \mathbf{x}^m = \sum\limits_{i_1, \cdots, i_m \in [n]} {a_{i_1 \cdots i_m} x_{i_1} \cdots x_{i_m}}$.
\end{lem}
\par
A great many methods have been proposed to compute $Z$-eigenpairs of symmetric tensors as yet. Qi et al. [\ref{Qi2}] developed a direct method to find all $Z$-eigenvalues of a third-order three dimension tensor. Kolda and Mayo [\ref{Kolda}] provided a shifted power method for computing all $Z$-eigenpairs, while Zeng and Ni [\ref{zeng}] proposed a quasi-Newton method to compute them. Cui et al. [\ref{cui}] considered Jacobian SDP relaxations in polynomial optimization for computing all $Z$-eigenpairs. Chen et al. [\ref{chen3}] proposed two homotopy continuation type algorithms to compute the $Z$-eigenpairs. Recently, dynamic method is also offered by Mo et al. [\ref{mo}] to reach this goal.

\par
Apart from the $Z$-eigenvalue, the $C$-eigenvalue of the piezoelectric-type tensor has attract many researchers' attention. Before going further, we recall the strict definitions of the above two terminologies in the following.
\begin{defn}\label{d1}
$([\ref{chen}])$ Let $\mathcal{A} = (a_{ijk} ) \in \mathbb{R}^{n  \times n  \times n}$ be a third-order $n$ dimensional real tensor. If the later two indices of $\mathcal{A}$ are symmetric, i.e., $a_{ijk} = a_{ikj}$ for all $j, k \in [n]$, then $\mathcal{A}$ is called a piezoelectric-type tensor.
\end{defn}
\par
Given a piezoelectric-type tensor $\mathcal{A}$, for each $i \in [n]$, its horizontal slice $\mathcal{A}(i, :, :) \in \mathbb{R}^{n \times n}$ denote the matrix whose $(j, k)$-entry is $a_{ijk}$. Then, it follows from Definition \ref{d1} that $\mathcal{A}(i, :, :)$ is a symmetric matrix.
\par
As we know, third order tensors have many important applications in some fields, such as physics and engineering [\ref{chen2}-\ref{Zou}]. As a special third-order tensor, piezoelectric-type tensors play a crucial role in piezoelectric effect and converse piezoelectric effect [\ref{chen}]. In order to investigate more properties related to piezoelectric effect and converse piezoelectric effect in solid crystal, the authors in [\ref{chen}] gave the definitions of $C$-eigenvalues and $C$-eigenvectors for piezoelectric-type tensors.

\begin{defn}\label{d2}
$([\ref{chen}])$ Let $\mathcal{A} = (a_{ijk} ) \in \mathbb{R}^{n  \times n  \times n}$ be a piezoelectric-type tensor. If there exist a scalar $\lambda \in \mathbb{R}$, and two vectors $\mathbf{x} =(x_1, x_2, \cdots, x_n)^T \in \mathbb{R}^n$ and $\mathbf{y} =(y_1, y_2, \cdots, y_n)^T \in \mathbb{R}^n$ satisfying the following system
$$\mathcal{A} \mathbf{y}\mathbf{y} = \lambda \mathbf{x}, \quad \mathbf{x} \mathcal{A} \mathbf{y} = \lambda \mathbf{y}, \quad \mathbf{x}^T \mathbf{x}=1, \quad \mathbf{y}^T \mathbf{y}=1,$$
where $\mathcal{A} \mathbf{y}\mathbf{y}  \in \mathbb{R}^n$ and $\mathbf{x} \mathcal{A} \mathbf{y}  \in \mathbb{R}^n$ with the $i$-th entry
$$(\mathcal{A} \mathbf{y}\mathbf{y})_i = \sum\limits_{j, k \in [n]} {a_{ijk} y_j y_k} \quad and \quad (\mathbf{x} \mathcal{A} \mathbf{y})_i = \sum\limits_{j, k \in [n]} {a_{jki} x_j y_k},$$
respectively, then $\lambda$ is called a $C$-eigenvalue of $\mathcal{A}$, and $\mathbf{x}$ and $\mathbf{y}$ are called associated left and right $C$-eigenvectors, respectively. Also, $(\lambda, \mathbf{x}, \mathbf{y})$ can be called a $C$-eigenpair of $\mathcal{A}$.
\end{defn}
\par
Suppose that $\mathbf{x}=(x_1, x_2, \cdots, x_n)^T$ and $\mathbf{y}=(y_1, y_2, \cdots, y_n)^T$ are the left and right $C$-eigenvectors corresponding to $\lambda$, then $\lambda = \mathbf{x} \mathcal{A} \mathbf{y}\mathbf{y} \triangleq \sum\limits_{i, j, k \in [n]} {a_{ijk} x_i y_j y_k}$. Moreover, $(\lambda, \mathbf{x}, -\mathbf{y})$, $(-\lambda, -\mathbf{x},  \mathbf{y})$, $(-\lambda, -\mathbf{x}, - \mathbf{y})$ are also $C$-eigenpairs of $\mathcal{A}$. For a given piezoelectric type tensor $\mathcal{A}$, let $\lambda_{C \mathop {\max }}(\mathcal{A})$ (or $\lambda_{C \mathop {\max }}$ for short) be the largest $C$-eigenvalue of $\mathcal{A}$, then
\begin{equation}\label{1.1}
\lambda_{C \mathop {\max }} = \max \{ \mathbf{x} \mathcal{A}  \mathbf{y} \mathbf{y}:\;  \mathbf{x}^T \mathbf{x}=1,  \mathbf{y}^T \mathbf{y}=1\}.
\end{equation}
Chen et al. [\ref{chen}] showed that the largest $C$-eigenvalue of the piezoelectric-type tensor has concrete physical meaning which determines the highest piezoelectric coupling constant.

\par
Due to the importance of the largest $C$-eigenvalue of the piezoelectric-type tensor, some authors focus their interest on this topic. As far as we know, there are several works related to compute all the $C$-eigenpairs of piezoelectric-type tensors. One is the work by Liang and Yang [\ref{liang}] in which they proposed a shifted eigenvalue decomposition method for computing $C$-eigenvalues. The others are proposed by Zhao and Luo [\ref{zhao}], Liu and Mo [\ref{liu2}] via $Z$-eigenpair methods. Yang and Liang [\ref{yang}] developed a convex relaxation method to compute the largest $C$-eigenvalue. Besides, some estimation are given  to locate all the $C$-eigenvalues in an interval, see [\ref{li}-\ref{liu1}] and the references therein.
\par

\par
Liu and Mo [\ref{liu2}] constructed a symmetric fourth-order tensor that related to a given piezoelectric-type tensor $\mathcal{A}$ by two steps. Firstly, construct a partially symmetric tensor $\mathcal{B}= (b_{i_1 i_2 i_3 i_4}) \in \mathbb{R}^{[4,n]}$ with entries defined as
\begin{equation}
b_{i_1 i_2 i_3 i_4} = \sum\limits_{i \in [n]} {a_{ii_1 i_2} a_{ii_3 i_4}},\;\; i_1, i_2, i_3,i_4 \in [n]. \label{1.02}
\end{equation}
In the next step, a symmetric tensor $\mathcal{\bar B} = (\bar b_{i_1 i_2 i_3 i_4})$ can be easily got from the partially symmetric tensor $\mathcal{B}$. We only need to set
\begin{equation}
\bar b_{i_1 i_2 i_3 i_4} = \frac{1}{3}(b_{i_1 i_2 i_3 i_4} + b_{i_1 i_3 i_2 i_4} + b_{i_1 i_4 i_2 i_3}). \label{1.2}
\end{equation}
They also showed the following close relationships between the $C$-eigenpair of the original piezoelectric-type tensor $\mathcal{A}$ and the $Z$-eigenpair of the structured tenor $\mathcal{\bar B}$.
\begin{lem}\label{l1}
$([\ref{liu2}])$ Let $\mathcal{A} = (a_{ijk} ) \in \mathbb{R}^{[3,n]}$ be a piezoelectric-type tensor, and $\mathcal{\bar B}$ be defined by $(\ref{1.2})$. \\
$(1)$ If $(\lambda, \mathbf{x}, \mathbf{y})$ is a $C$-eigenpair of $\mathcal{A}$, then $(\lambda^2, \mathbf{y})$ is a $Z$-eigenpair of $\mathcal{\bar B}$.\\
$(2)$ If $(\lambda^2, \mathbf{y})$ is a $Z$-eigenpair of $\mathcal{\bar B}$ with $\lambda \ne 0$, then $(\lambda, \mathbf{x}, \mathbf{y})$ is a $C$-eigenpair of $\mathcal{A}$ with $\mathbf{x} = \frac{1}{\lambda}\mathcal{A} \mathbf{y}\mathbf{y}$. \\
$(3)$ If $(0, \mathbf{y})$ is a $Z$-eigenpair of $\mathcal{\bar B}$, then $(0, \mathbf{x}, \mathbf{y})$ is a $C$-eigenpair of $\mathcal{A}$, where $\mathbf{x}$ is a real unit solution of the equation $\mathbf{x} \mathcal{A}  \mathbf{y} = 0$.
\end{lem}
{
\par
Furthermore, we can prove the following result.
\begin{thm}\label{th1.1}
Let $\mathcal{A} = (a_{ijk} ) \in \mathbb{R}^{[3,n]}$ be a piezoelectric-type tensor, $\mathcal{B}$ and $\mathcal{\bar B}$ be defined by $(\ref{1.02})$ and $(\ref{1.2})$. Then $\mathcal{\bar B}$ is a positive semi-definite $(PSD)$ tensor, i.e.,
$$\mathcal{\bar B} \mathbf{x}^4 = \sum\limits_{i_1, i_2, i_3, i_4 \in [n]} {\bar b_{i_1 i_2 i_3 i_4} x_{i_1} x_{i_2} x_{i_3} x_{i_4}} \geqslant 0 \;\;for\;all\; \mathbf{x} \in \mathbb{R}^{n}.$$
Therefore, all the $Z$-eigenvalues of $\mathcal{\bar B}$ are nonnegative.
\end{thm}
\textbf{Proof.} For each $\mathbf{x} \in \mathbb{R}^{n}$, Liu and Mo [\ref{liu2}] proved that $\mathcal{\bar B} \mathbf{x}^4 = \mathcal{B} \mathbf{x}^4$. Hence,
\begin{eqnarray*}
\mathcal{\bar B} \mathbf{x}^4 &=& \mathcal{B} \mathbf{x}^4 = \sum\limits_{i_1, i_2, i_3, i_4 \in [n]} {b_{i_1 i_2 i_3 i_4} x_{i_1} x_{i_2} x_{i_3} x_{i_4}} \\
&=& \sum\limits_{i_1, i_2, i_3, i_4 \in [n]} {\left( {\sum\limits_{i \in [n]} {a_{ii_1 i_2} a_{ii_3 i_4}}} \right)} x_{i_1} x_{i_2} x_{i_3} x_{i_4} \\
&=& \sum\limits_{i \in [n]} {(\mathcal{A} \mathbf{x} \mathbf{x})_i^2} \geqslant 0.
\end{eqnarray*}
It follows from [Theorem 5, \ref{Qi1}] that all the $Z$-eigenvalues of $\mathcal{\bar B}$ are nonnegative. The proof is completed. \;\;\;$\Box$
}

\par
Actually, $\mathcal{\bar B}$ can be seen as a function of piezoelectric-type tensor $\mathcal{A}$. For convenience, for a given piezoelectric-type tensor $\mathcal{A}$, we adopt $\mathcal{S}_{\mathcal{A}}$ to denote this function, namely, $\mathcal{S}_{\mathcal{A}} = \mathcal{\bar B}$.

\section{Main results}

Let $\mathcal{A}$ and its perturbed tensor $ \mathcal{\tilde A} = \mathcal{A} + \mathcal{E}$ be piezoelectric-type tensors. According to the definition of piezoelectric-type tensors, we know that the tensor $\mathcal{E}$ is also a piezoelectric-type tensor. Based on (\ref{1.1}), we have the following perturbation bound.

\begin{thm}\label{th2.1}
Let $\mathcal{A}$ and its perturbed tensor $ \mathcal{\tilde A} = \mathcal{A} + \mathcal{E}$ be piezoelectric-type tensors. Then
\begin{equation}
\lambda_{C \mathop {\max }} (\mathcal{A}) - \lambda_{C \mathop {\max }} (\mathcal{E}) \leqslant \lambda_{C \mathop {\max }} (\mathcal{\tilde A}) \leqslant \lambda_{C \mathop {\max }} (\mathcal{A}) + \lambda_{C \mathop {\max }} (\mathcal{E}). \label{th2.1-01}
\end{equation}
\end{thm}
\textbf{Proof.} In view of (\ref{1.1}), we have
\begin{eqnarray}
	\lambda_{C \mathop {\max }} (\mathcal{\tilde A}) &=& \max \{ \mathbf{x} \mathcal{\tilde A} \mathbf{y}\mathbf{y}:\; \mathbf{x}^T \mathbf{x} = 1,\; \mathbf{y}^T \mathbf{y} = 1\} \nonumber\\
	&=& \max \{ \mathbf{x}\mathcal{A} \mathbf{y}\mathbf{y} + \mathbf{x} \mathcal{E} \mathbf{y}\mathbf{y}:\; \mathbf{x}^T \mathbf{x} = 1,\; \mathbf{y}^T \mathbf{y} = 1\} \nonumber\\
	&\leqslant& \max \{ \mathbf{x} \mathcal{A} \mathbf{y}\mathbf{y} :\; \mathbf{x}^T \mathbf{x}= 1,\; \mathbf{y}^T \mathbf{y} = 1\} + \max \{\mathbf{x} \mathcal{E} \mathbf{y}\mathbf{y}:\; \mathbf{x}^T \mathbf{x} = 1,\; \mathbf{y}^T \mathbf{y} = 1\} \nonumber\\
	&=& \lambda_{C \mathop {\max }}(\mathcal{A}) + \lambda_{C \mathop {\max }} (\mathcal{E}). \label{th2.1-1}
\end{eqnarray}
On the other hand, let $\mathbf{x}$ and $\mathbf{y}$ be the associated left and right $C$-eigenvectors corresponding to $\lambda_{C \mathop {\max }} (\mathcal{A})$ of piezoelectric-type $\mathcal{A}$. Then
\begin{equation}
	\lambda_{C \mathop {\max }} (\mathcal{\tilde A}) \geqslant \mathbf{x} \mathcal{\tilde A} \mathbf{y}\mathbf{y} = \mathbf{x} \mathcal{A} \mathbf{y}\mathbf{y} + \mathbf{x} \mathcal{E} \mathbf{y}\mathbf{y}= \lambda_{C \mathop {\max }} (\mathcal{A}) + \mathbf{x} \mathcal{E} \mathbf{y}\mathbf{y} \geqslant \lambda_{C \mathop {\max }} (\mathcal{A}) - \lambda_{C \mathop {\max }} (\mathcal{E}). \label{th2.1-2}
\end{equation}
Combining (\ref{th2.1-1}) and (\ref{th2.1-2}), we get the perturbation bound (\ref{th2.1-01}). \;\;\;$\Box$

\par
In [\ref{liu1}], the authors gave an effective estimation (upper bound) on the lagest $C$-eigenvalue of piezoelectric-type tensor. Hence, together with Theorem \ref{th2.1} and [Theorem 2.4, \ref{liu1}] produces the following result immediately.

\begin{cor}\label{c2.1}
Let $\mathcal{A}$ and its perturbed tensor $ \mathcal{\tilde A} = \mathcal{A} + \mathcal{E}$ be piezoelectric-type tensors. Then
\begin{equation}
\lambda_{C \mathop {\max }} (\mathcal{A}) - \left\| E \right\|_2 \leqslant \lambda_{C \mathop {\max }} (\mathcal{\tilde A}) \leqslant \lambda_{C \mathop {\max }} (\mathcal{A}) + \left\| E \right\|_2, \label{c2.1-01}
\end{equation}
where $E = \left( {\begin{array}{*{20}c}
   {\mathcal{E}(1, :, :)}, &  \cdots,  & {\mathcal{E}(n, :, :)} \\
 \end{array} } \right)$ is a column block matrix with the $i$-th block is $\mathcal{E}(i, :, :)$.
\end{cor}
\par
By the result [Theorem 2.4, \ref{liu1}], we know $\lambda_{C \mathop {\max }} (\mathcal{E}) \leqslant \left\| E \right\|_2$. Thus, we have
$$[\lambda_{C \mathop {\max }} (\mathcal{A}) - \lambda_{C \mathop {\max }} (\mathcal{E}),\; \lambda_{C \mathop {\max }} (\mathcal{A}) + \lambda_{C \mathop {\max }} (\mathcal{E})] \subseteq [\lambda_{C \mathop {\max }} (\mathcal{A}) - \left\| E \right\|_2,\;\lambda_{C \mathop {\max }} (\mathcal{A}) + \left\| E \right\|_2].$$

\par
By the relationships between the original piezoelectric-type tensor $\mathcal{A}$ and the structured tenor $\mathcal{S}_{\mathcal{A}}$ (Lemma \ref{l1}), we can deduce the following alternative perturbation bound.

\begin{thm}\label{th2.2}
Let $\mathcal{A}$ and its perturbed tensor $ \mathcal{\tilde A} = \mathcal{A} + \mathcal{E}$ be piezoelectric-type tensors. Then
\begin{equation}
\sqrt {(\lambda_{C \mathop {\max }} (\mathcal{A}))^2 + \lambda_{Z \mathop {\min }} (\mathcal{S}_{\mathcal{\tilde A}} - \mathcal{S}_{\mathcal{\mathcal{A}}})} \leqslant \lambda_{C \mathop {\max }} (\mathcal{\tilde A}) \leqslant \sqrt {(\lambda_{C \mathop {\max }} (\mathcal{A}))^2 + \lambda_{Z \mathop {\max }} (\mathcal{S}_{\mathcal{\tilde A}} - \mathcal{S}_{\mathcal{\mathcal{A}}})}. \label{th2.2-01}
\end{equation}
\end{thm}
\textbf{Proof.} In view of Lemma \ref{l1}, we have
\begin{eqnarray}
	(\lambda_{C \mathop {\max }} (\mathcal{\tilde A}))^2 &=& \lambda_{Z \mathop {\max }} (\mathcal{S}_{\mathcal{\tilde A}}) = \max \{ (\mathcal{S}_{\mathcal{\tilde A}}) \mathbf{y}^4:\;  \mathbf{y}^T \mathbf{y} = 1\} \nonumber\\
	&=& \max \{ (\mathcal{S}_{\mathcal{\mathcal{A}}}) \mathbf{y}^4 + (\mathcal{S}_{\mathcal{\tilde A}} - \mathcal{S}_{\mathcal{\mathcal{A}}}) \mathbf{y}^4:\;  \mathbf{y}^T \mathbf{y} = 1\} \nonumber\\
	&\leqslant& \max \{ (\mathcal{S}_{\mathcal{\mathcal{A}}}) \mathbf{y}^4:\;  \mathbf{y}^T \mathbf{y} = 1\} + \max \{(\mathcal{S}_{\mathcal{\tilde A}} - \mathcal{S}_{\mathcal{\mathcal{A}}}) \mathbf{y}^4:\;  \mathbf{y}^T \mathbf{y} = 1\} \nonumber\\
	&=& \lambda_{Z \mathop {\max }} (\mathcal{S}_{\mathcal{A}}) + \lambda_{Z \mathop {\max }} (\mathcal{S}_{\mathcal{\tilde A}} - \mathcal{S}_{\mathcal{\mathcal{A}}}) \nonumber\\
	&=& (\lambda_{C \mathop {\max }} (\mathcal{A}))^2 + \lambda_{Z \mathop {\max }} (\mathcal{S}_{\mathcal{\tilde A}} - \mathcal{S}_{\mathcal{\mathcal{A}}}). \label{th2.2-1}
\end{eqnarray}
On the other hand, let $\mathbf{y}_0$ be the associated right $C$-eigenvector corresponding to $\lambda_{C \mathop {\max }} (\mathcal{A})$ of piezoelectric-type $\mathcal{A}$. Clearly, $\mathbf{y}_0$ is also the associated $Z$-eigenvector corresponding to $\lambda_{Z \mathop {\max }} (\mathcal{S}_{\mathcal{A}})$ of $\mathcal{S}_{\mathcal{A}}$. Then

\begin{eqnarray}
	(\lambda_{C \mathop {\max }} (\mathcal{\tilde A}))^2 &=& \lambda_{Z \mathop {\max }} (\mathcal{S}_{\mathcal{\tilde A}}) = \max \{ (\mathcal{S}_{\mathcal{\tilde A}}) \mathbf{y}^4:\;  \mathbf{y}^T \mathbf{y} = 1\} \nonumber\\
	&\geqslant& (\mathcal{S}_{\mathcal{\mathcal{A}}}) \mathbf{y}_0^4 + (\mathcal{S}_{\mathcal{\tilde A}} - \mathcal{S}_{\mathcal{\mathcal{A}}}) \mathbf{y}_0^4 \nonumber\\
    &=& \lambda_{Z \mathop {\max }} (\mathcal{S}_{\mathcal{A}}) + (\mathcal{S}_{\mathcal{\tilde A}} - \mathcal{S}_{\mathcal{\mathcal{A}}}) \mathbf{y}_0^4  \nonumber\\
	&\geqslant& \lambda_{Z \mathop {\max }} (\mathcal{S}_{\mathcal{A}}) + \lambda_{Z \mathop {\min}} (\mathcal{S}_{\mathcal{\tilde A}} - \mathcal{S}_{\mathcal{\mathcal{A}}})  \nonumber\\
	&=& (\lambda_{C \mathop {\max }} (\mathcal{A}))^2 + \lambda_{Z \mathop {\min}} (\mathcal{S}_{\mathcal{\tilde A}} - \mathcal{S}_{\mathcal{\mathcal{A}}}). \label{th2.2-2}
\end{eqnarray}
{
Moreover, for each $\mathbf{y}$ with $\mathbf{y}^T \mathbf{y} = 1$, we have
$$(\mathcal{S}_{\mathcal{\tilde A}} - \mathcal{S}_{\mathcal{A}}) \mathbf{y}^4 = \mathcal{S}_{\mathcal{\tilde A}} \mathbf{y}^4 - \mathcal{S}_{\mathcal{A}} \mathbf{y}^4 \geqslant \mathcal{S}_{\mathcal{\tilde A}} \mathbf{y}^4 - \lambda_{Z \mathop {\max }} (\mathcal{S}_{\mathcal{A}}) = \mathcal{S}_{\mathcal{\tilde A}} \mathbf{y}^4 - (\lambda_{C \mathop {\max }} (\mathcal{A}))^2.$$
Applying Theorem \ref{th1.1} to the above inequality, we have
$$(\mathcal{S}_{\mathcal{\tilde A}} - \mathcal{S}_{\mathcal{A}}) \mathbf{y}^4 + (\lambda_{C \mathop {\max }} (\mathcal{A}))^2 \geqslant \mathcal{S}_{\mathcal{\tilde A}} \mathbf{y}^4 \geqslant 0,$$
which means that both $(\lambda_{C \mathop {\max }} (\mathcal{A}))^2 + \lambda_{Z \mathop {\max }} (\mathcal{S}_{\mathcal{\tilde A}} - \mathcal{S}_{\mathcal{\mathcal{A}}})$ and $(\lambda_{C \mathop {\max }} (\mathcal{A}))^2 + \lambda_{Z \mathop {\min }} (\mathcal{S}_{\mathcal{\tilde A}} - \mathcal{S}_{\mathcal{\mathcal{A}}})$ are nonnegative.
}
\par
Combining (\ref{th2.2-1}) and (\ref{th2.2-2}), we get the perturbation bound (\ref{th2.2-01}). \;\;\;$\Box$

{
\par
The following theorem shows that the perturbation bound given by Theorem \ref{th2.2} is more accurate than the one given by Theorem \ref{th2.1}.
\begin{thm}\label{th2.3}
Let $\mathcal{A}$ and its perturbed tensor $ \mathcal{\tilde A} = \mathcal{A} + \mathcal{E}$ be piezoelectric-type tensors. Then
\begin{eqnarray}
&&\left[ {\sqrt {(\lambda_{C \mathop {\max }} (\mathcal{A}))^2 + \lambda_{Z \mathop {\min }} (\mathcal{S}_{\mathcal{\tilde A}} - \mathcal{S}_{\mathcal{\mathcal{A}}})},\;\sqrt {(\lambda_{C \mathop {\max }} (\mathcal{A}))^2 + \lambda_{Z \mathop {\max }} (\mathcal{S}_{\mathcal{\tilde A}} - \mathcal{S}_{\mathcal{\mathcal{A}}})}}\; \right] \nonumber \\
&\subseteq& [\lambda_{C \mathop {\max }} (\mathcal{A}) - \lambda_{C \mathop {\max }} (\mathcal{E}),\; \lambda_{C \mathop {\max }} (\mathcal{A}) + \lambda_{C \mathop {\max }} (\mathcal{E})]. \label{th2.3-01}
\end{eqnarray}
\end{thm}
\textbf{Proof.} For each piezoelectric-type tensor $\mathcal{A}$, it follows from [Theorem 2.1, \ref{liu2}] that the fourth-order tensors $\mathcal{B}$ and $\mathcal{\bar B}$ defined by (\ref{1.02}) and (\ref{1.2}) have the same $Z$-eigenpairs. So, without loss of generality, we assume that $\mathcal{S}_{\mathcal{A}} = \mathcal{B}$.
\par
Firstly, we come to prove
$$\sqrt {(\lambda_{C \mathop {\max }} (\mathcal{A}))^2 + \lambda_{Z \mathop {\max }} (\mathcal{S}_{\mathcal{\tilde A}} - \mathcal{S}_{\mathcal{\mathcal{A}}})} \leqslant \lambda_{C \mathop {\max }} (\mathcal{A}) + \lambda_{C \mathop {\max }} (\mathcal{E}).$$
In view of Lemma \ref{l1}, the above inequality is equivalent to
\begin{eqnarray}
\lambda_{Z \mathop {\max }} (\mathcal{S}_{\mathcal{\tilde A}} - \mathcal{S}_{\mathcal{\mathcal{A}}}) &\leqslant& (\lambda_{C \mathop {\max }} (\mathcal{E}))^2 + 2 \lambda_{C \mathop {\max }} (\mathcal{A})\lambda_{C \mathop {\max }} (\mathcal{E}) \nonumber \\
&=& \lambda_{Z \mathop {\max }} (\mathcal{S}_{\mathcal{E}}) + 2 \lambda_{C \mathop {\max }} (\mathcal{A})\lambda_{C \mathop {\max }} (\mathcal{E}). \label{th2.3-1}
\end{eqnarray}
Let $\mathcal{A} = (a_{ijk} ) \in \mathbb{R}^{[3,n]}$ and $\mathcal{E} = (\varepsilon_{ijk} ) \in \mathbb{R}^{[3,n]}$, then for each vector $\mathbf{y}$ with $\mathbf{y}^T \mathbf{y} = 1$, we have
\begin{eqnarray*}
&& (\mathcal{S}_{\mathcal{\tilde A}} - \mathcal{S}_{\mathcal{A}}) \mathbf{y}^4 \\
&=& \sum\limits_{i_1, i_2, i_3, i_4 \in [n]} {\left( {\sum\limits_{i \in [n]} {(a_{ii_1 i_2} + \varepsilon_{ii_1 i_2})(a_{ii_3 i_4} + \varepsilon_{ii_3 i_4})}  - \sum\limits_{i \in [n]} {a_{ii_1 i_2} a_{ii_3 i_4}} } \right)} y_{i_1} y_{i_2} y_{i_3} y_{i_4}\\
&=& \sum\limits_{i_1, i_2, i_3, i_4 \in [n]} {\left( {\sum\limits_{i \in [n]} {\varepsilon_{ii_1 i_2} \varepsilon_{ii_3 i_4}}} \right)} y_{i_1} y_{i_2} y_{i_3} y_{i_4} + 2 \sum\limits_{i_1, i_2, i_3, i_4 \in [n]} {\left( {\sum\limits_{i \in [n]} {a_{ii_1 i_2} \varepsilon_{ii_3 i_4}}} \right)} y_{i_1} y_{i_2} y_{i_3} y_{i_4} \\
&=& \mathcal{S}_{\mathcal{E}} \mathbf{y}^4 + 2 \sum\limits_{i \in [n]} {(\mathcal{A} \mathbf{y} \mathbf{y})_i (\mathcal{E} \mathbf{y} \mathbf{y})_i}  \\
&\leqslant& \mathcal{S}_{\mathcal{E}} \mathbf{y}^4  + 2 \sqrt {\sum\limits_{i \in [n]} {(\mathcal{A} \mathbf{y} \mathbf{y})_i^2 }}  \sqrt {\sum\limits_{i \in [n]} {(\mathcal{E} \mathbf{y} \mathbf{y})_i^2}} \\
&=& \mathcal{S}_{\mathcal{E}} \mathbf{y}^4 + 2 \sqrt {\mathcal{S}_{\mathcal{A}} \mathbf{y}^4}  \sqrt {\mathcal{S}_{\mathcal{E}} \mathbf{y}^4}.
\end{eqnarray*}
Therefore,
\begin{eqnarray*}
&& \max \{ (\mathcal{S}_{\mathcal{\tilde A}} - \mathcal{S}_{\mathcal{A}}) \mathbf{y}^4:\;  \mathbf{y}^T \mathbf{y} = 1\} \leqslant  \max \{ \mathcal{S}_{\mathcal{E}} \mathbf{y}^4 + 2 \sqrt {\mathcal{S}_{\mathcal{A}} \mathbf{y}^4}  \sqrt {\mathcal{S}_{\mathcal{E}} \mathbf{y}^4}:\;  \mathbf{y}^T \mathbf{y} = 1\} \\
&\leqslant& \max \{\mathcal{S}_{\mathcal{E}} \mathbf{y}^4:\;  \mathbf{y}^T \mathbf{y} = 1\} + \max \{2 \sqrt {\mathcal{S}_{\mathcal{A}} \mathbf{y}^4}  \sqrt {\mathcal{S}_{\mathcal{E}} \mathbf{y}^4}:\;  \mathbf{y}^T \mathbf{y} = 1\},
\end{eqnarray*}
which produces
$$\lambda_{Z \mathop {\max }} (\mathcal{S}_{\mathcal{\tilde A}} - \mathcal{S}_{\mathcal{A}}) \leqslant \lambda_{Z \mathop {\max }} (\mathcal{S}_{\mathcal{E}}) + 2\sqrt {\lambda_{Z \mathop {\max }} (\mathcal{S}_\mathcal{A})} \sqrt {\lambda_{Z \mathop {\max }} (\mathcal{S}_\mathcal{E})}.$$
Applying Lemma \ref{l1} to the above inequality, we get (\ref{th2.3-1}).
\par
Next, we prove
\begin{equation}
\lambda_{C \mathop {\max }} (\mathcal{A}) - \lambda_{C \mathop {\max }} (\mathcal{E}) \leqslant \sqrt {(\lambda_{C \mathop {\max }} (\mathcal{A}))^2 + \lambda_{Z \mathop {\min }} (\mathcal{S}_{\mathcal{\tilde A}} - \mathcal{S}_{\mathcal{\mathcal{A}}})}.\label{th2.3-02}
\end{equation}
When $\lambda_{C \mathop {\max }} (\mathcal{A}) < \lambda_{C \mathop {\max }} (\mathcal{E})$, (\ref{th2.3-02}) is evident. When $\lambda_{C \mathop {\max }} (\mathcal{A}) \geqslant \lambda_{C \mathop {\max }} (\mathcal{E})$, (\ref{th2.3-02}) is equivalent to
\begin{eqnarray}
\lambda_{Z \mathop {\min }} (\mathcal{S}_{\mathcal{\tilde A}} - \mathcal{S}_{\mathcal{\mathcal{A}}}) &\geqslant& (\lambda_{C \mathop {\max }} (\mathcal{E}))^2 - 2 \lambda_{C \mathop {\max }} (\mathcal{A})\lambda_{C \mathop {\max }} (\mathcal{E}) \nonumber \\
&=& \lambda_{Z \mathop {\max }} (\mathcal{S}_{\mathcal{E}}) - 2 \lambda_{C \mathop {\max }} (\mathcal{A})\lambda_{C \mathop {\max }} (\mathcal{E}). \label{th2.3-2}
\end{eqnarray}
Since
\begin{eqnarray*}
(\mathcal{S}_{\mathcal{\tilde A}} - \mathcal{S}_{\mathcal{A}}) \mathbf{y}^4 &=& \mathcal{S}_{\mathcal{E}} \mathbf{y}^4 + 2 \sum\limits_{i \in [n]} {(\mathcal{A} \mathbf{y} \mathbf{y})_i (\mathcal{E} \mathbf{y} \mathbf{y})_i} \\
&\geqslant& \mathcal{S}_{\mathcal{E}} \mathbf{y}^4 - 2 \sqrt {\sum\limits_{i \in [n]} {(\mathcal{A} \mathbf{y} \mathbf{y})_i^2 }}  \sqrt {\sum\limits_{i \in [n]} {(\mathcal{E} \mathbf{y} \mathbf{y})_i^2}} \\
&=& \mathcal{S}_{\mathcal{E}} \mathbf{y}^4 - 2 \sqrt {\mathcal{S}_{\mathcal{A}} \mathbf{y}^4}  \sqrt {\mathcal{S}_{\mathcal{E}} \mathbf{y}^4}.
\end{eqnarray*}
Apply Lemma \ref{l1} again, and consider the following difference,
\begin{eqnarray*}
&& \left( {\mathcal{S}_{\mathcal{E}} \mathbf{y}^4 - 2 \sqrt {\mathcal{S}_{\mathcal{A}} \mathbf{y}^4}  \sqrt {\mathcal{S}_{\mathcal{E}} \mathbf{y}^4}} \right) - [\lambda_{Z \mathop {\max }} (\mathcal{S}_{\mathcal{E}}) - 2 \lambda_{C \mathop {\max }} (\mathcal{A})\lambda_{C \mathop {\max }} (\mathcal{E})] \nonumber\\
&=& [\mathcal{S}_{\mathcal{E}} \mathbf{y}^4 - \lambda_{Z \mathop {\max }} (\mathcal{S}_{\mathcal{E}})] +2 \left( {\sqrt {\lambda_{Z \mathop {\max }} (\mathcal{S}_\mathcal{A})} \sqrt {\lambda_{Z \mathop {\max }} (\mathcal{S}_\mathcal{E})} - \sqrt {\mathcal{S}_{\mathcal{A}} \mathbf{y}^4}  \sqrt {\mathcal{S}_{\mathcal{E}} \mathbf{y}^4}} \right) \nonumber\\
&=& [\mathcal{S}_{\mathcal{E}} \mathbf{y}^4 - \lambda_{Z \mathop {\max }} (\mathcal{S}_{\mathcal{E}})] - 2 \sqrt {\lambda_{Z \mathop {\max }} (\mathcal{S}_\mathcal{A})} \left( {\sqrt {\mathcal{S}_{\mathcal{E}} \mathbf{y}^4} - \sqrt {\lambda_{Z \mathop {\max }} (\mathcal{S}_\mathcal{E})}} \right) \nonumber\\
&& + 2 \sqrt {\mathcal{S}_{\mathcal{E}} \mathbf{y}^4} \left( {\sqrt {\lambda_{Z \mathop {\max }} (\mathcal{S}_\mathcal{A})} - \sqrt {\mathcal{S}_{\mathcal{A}} \mathbf{y}^4}} \right) \nonumber\\
&=& \left( {\sqrt {\mathcal{S}_{\mathcal{E}} \mathbf{y}^4} - \sqrt {\lambda_{Z \mathop {\max }} (\mathcal{S}_\mathcal{E})}} \right) \left( {\sqrt {\mathcal{S}_{\mathcal{E}} \mathbf{y}^4} + \sqrt {\lambda_{Z \mathop {\max }} (\mathcal{S}_\mathcal{E})} - 2 \sqrt {\lambda_{Z \mathop {\max }} (\mathcal{S}_\mathcal{A})}} \right)   \nonumber\\
&& + 2 \sqrt {\mathcal{S}_{\mathcal{E}} \mathbf{y}^4} \left( {\sqrt {\lambda_{Z \mathop {\max }} (\mathcal{S}_\mathcal{A})} - \sqrt {\mathcal{S}_{\mathcal{A}} \mathbf{y}^4}} \right).
\end{eqnarray*}
Note that $\lambda_{C \mathop {\max }} (\mathcal{A}) \geqslant \lambda_{C \mathop {\max }} (\mathcal{E})$, $\mathcal{S}_{\mathcal{E}} \mathbf{y}^4 \leqslant \lambda_{Z \mathop {\max }} (\mathcal{S}_\mathcal{E})$ and $\mathcal{S}_{\mathcal{A}} \mathbf{y}^4 \leqslant \lambda_{Z \mathop {\max }} (\mathcal{S}_\mathcal{A})$, then, we have
$$\left( {\mathcal{S}_{\mathcal{E}} \mathbf{y}^4 - 2 \sqrt {\mathcal{S}_{\mathcal{A}} \mathbf{y}^4}  \sqrt {\mathcal{S}_{\mathcal{E}} \mathbf{y}^4}} \right) - [\lambda_{Z \mathop {\max }} (\mathcal{S}_{\mathcal{E}}) - 2 \lambda_{C \mathop {\max }} (\mathcal{A})\lambda_{C \mathop {\max }} (\mathcal{E})] \geqslant 0,$$
i.e.,
$$\mathcal{S}_{\mathcal{E}} \mathbf{y}^4 - 2 \sqrt {\mathcal{S}_{\mathcal{A}} \mathbf{y}^4}  \sqrt {\mathcal{S}_{\mathcal{E}} \mathbf{y}^4} \geqslant \lambda_{Z \mathop {\max }} (\mathcal{S}_{\mathcal{E}}) - 2 \lambda_{C \mathop {\max }} (\mathcal{A})\lambda_{C \mathop {\max }} (\mathcal{E}).$$
Therefore,
$$(\mathcal{S}_{\mathcal{ A}} - \mathcal{S}_{\mathcal{\tilde A}}) \mathbf{y}^4 \geqslant \lambda_{Z \mathop {\max }} (\mathcal{S}_{\mathcal{E}}) - 2 \lambda_{C \mathop {\max }} (\mathcal{A})\lambda_{C \mathop {\max }} (\mathcal{E}),$$
which produces (\ref{th2.3-2}) and (\ref{th2.3-02}).
\par
By the above analyses, we get the desired result (\ref{th2.3-01}).   \;\;\;$\Box$

\par
Comparing with the bounds in Theorem \ref{th2.1}, Corollary \ref{c2.1} and Theorem \ref{th2.2}, we have proved that the perturbation bound in Theorem \ref{th2.2} is always tightest than the other ones.
}

\section{Numerical examples}
In this section, we will give some numerical examples to illustrate the efficiency of the proposed bounds. The tested examples are from the real world, detailedly, we can see Examples 1-8 in \cite{chen}. Let $\mathcal{E}$ generated by the following codes of Matlab.
\begin{equation*}
	\hat{\mathcal{E}}=\epsilon\textup{rand}(n,n,n),~\mathcal{E}=\textup{symmetrize}(\textup{tensor}(\hat{\mathcal{E}}),[2~ 3]),
\end{equation*}
where \textit{symmetrize} and  \textit{tensor} are the functions of 'Tensor Toolbox 2.6' (see the work of Bader et al. [\ref{ma1}]). The true $C$-eigenvalue of $\mathcal{\tilde A}$ is computed by the \textit{zeig} function of the toolbox 'TenEig', which is a MATLAB toolbox to find eigenpairs of a tensor given by the work of [\ref{ma2}].

\begin{table}[htbp]
	\centering
  \caption{Comparison the proposed upper bounds with the true $C$-eigenvalue}
	\scalebox{0.58}{
	\begin{tabular}{lcccccccc}
\hline
		$\epsilon$     &       & 1     & $10^{-1}$   & $10^{-2}$  & $10^{-3}$ & $10^{-4}$ & $10^{-5}$ \\
		\hline
    $VFeSb$ & TRUE  & 4.87650627  & 4.30866681  & 4.25706047  & 4.25195125  & 4.25144085  & 4.25138981  \\
          & (2.1) & 7.14734104  & 4.54097983  & 4.28034371  & 4.25428010  & 4.25167374  & 4.25141310  \\
          & (2.4) & 7.16334615  & 4.54258034  & 4.28050376  & 4.25429611  & 4.25167534  & 4.25141326  \\
          & (2.5) & \textbf{4.89007959 } & \textbf{4.30878542 } & \textbf{4.25706736 } & \textbf{4.25195190 } & \textbf{4.25144091 } & \textbf{4.25138982 } \\
          \hline
    $SiO_2$  & TRUE  & 2.98257552  & 0.39023942  & 0.15777358  & 0.13896998  & 0.13762874  & 0.13754487  \\
          & (2.1) & 3.03349312  & 0.42713192  & 0.16649580  & 0.14043218  & 0.13782582  & 0.13756519  \\
          & (2.4) & 3.04949823  & 0.42873243  & 0.16665585  & 0.14044819  & 0.13782742  & 0.13756535  \\
          & (2.5) & \textbf{2.98359085 } & \textbf{0.39521633 } & \textbf{0.15999258 } & \textbf{0.13971317 } & \textbf{0.13775317 } & \textbf{0.13755791 } \\
          \hline
    $Cr_2AgBiO_8$ & TRUE  & 3.74259460  & 2.70447034  & 2.63169114  & 2.62632963  & 2.62585024  & 2.62580287  \\
          & (2.1) & 5.52175451  & 2.91539330  & 2.65475718  & 2.62869357  & 2.62608721  & 2.62582657  \\
          & (2.4) & 5.53775962  & 2.91699381  & 2.65491723  & 2.62870957  & 2.62608881  & 2.62582673  \\
          & (2.5) & \textbf{3.77276365 } & \textbf{2.70944409 } & \textbf{2.63393422 } & \textbf{2.62660903 } & \textbf{2.62587873 } & \textbf{2.62580572 } \\
          \hline
    $RbTaO_3$ & TRUE  & 13.01446669  & 13.56705411  & 13.63089784  & 13.63738109  & 13.63803042  & 13.63809536  \\
          & (2.1) & 16.53405947  & 13.92769827  & 13.66706215  & 13.64099853  & 13.63839217  & 13.63813154  \\
          & (2.4) & 16.55006458  & 13.92929878  & 13.66722220  & 13.64101454  & 13.63839377  & 13.63813170  \\
          & (2.5) & \textbf{14.53398426 } & \textbf{13.72158318 } & \textbf{13.64639155 } & \textbf{13.63893089 } & \textbf{13.63818540 } & \textbf{13.63811086 } \\
          \hline
    $NaBiS_2$ & TRUE  & 11.57301393  & 11.65711348  & 11.66633731  & 11.66726751  & 11.66736061  & 11.66736992  \\
          & (2.1) & 14.56332785  & 11.95696665  & 11.69633052  & 11.67026691  & 11.66766055  & 11.66739992  \\
          & (2.4) & 14.57933296  & 11.95856716  & 11.69649058  & 11.67028292  & 11.66766215  & 11.66740008  \\
          & (2.5) & \textbf{11.95961970 } & \textbf{11.69474089 } & \textbf{11.67008916 } & \textbf{11.66764259 } & \textbf{11.66739812 } & \textbf{11.66737367 } \\
          \hline
    $LiBiB_2O_5$ & TRUE  & 9.48563900  & 7.88963674  & 7.75257003  & 7.73911060  & 7.73776715  & 7.73763283  \\
          & (2.1) & 10.63357480  & 8.02721360  & 7.76657748  & 7.74051386  & 7.73790750  & 7.73764687  \\
          & (2.4) & 10.64957991  & 8.02881411  & 7.76673753  & 7.74052987  & 7.73790910  & 7.73764703  \\
          & (2.5) & \textbf{9.59517350 } & \textbf{7.89942442 } & \textbf{7.75351783 } & \textbf{7.73920505 } & \textbf{7.73777659 } & \textbf{7.73763378 } \\
          \hline
    $KBi_2F_7$ & TRUE  & 14.68094239  & 13.60467138  & 13.51224483  & 13.50314448  & 13.50223585  & 13.50214501  \\
          & (2.1) & 16.39809181  & 13.79173060  & 13.53109448  & 13.50503087  & 13.50242451  & 13.50216387  \\
          & (2.4) & 16.41409692  & 13.79333111  & 13.53125453  & 13.50504687  & 13.50242611  & 13.50216403  \\
          & (2.5) & \textbf{15.26791582 } & \textbf{13.66616563 } & \textbf{13.51840457 } & \textbf{13.50376054 } & \textbf{13.50229746 } & \textbf{13.50215117 } \\
          \hline
    $BaNiO_3$ & TRUE  & 27.54177747  & 27.46666969  & 27.46314869  & 27.46283449  & 27.46280345  & 27.46280034  \\
          & (2.1) & 30.35875689  & 27.75239569  & 27.49175957  & 27.46569596  & 27.46308960  & 27.46282896  \\
          & (2.4) & 30.37476201  & 27.75399620  & 27.49191962  & 27.46571196  & 27.46309120  & 27.46282912  \\
          & (2.5) & \textbf{28.28527171 } & \textbf{27.53589622 } & \textbf{27.47002036 } & \textbf{27.46352115 } & \textbf{27.46287211 } & \textbf{27.46280721 } \\
    \hline
    \end{tabular}}%
  \label{tab_up}%
\end{table}%
\begin{table}[htbp]
  \centering
	\caption{Comparison the proposed lower bounds with the true $C$-eigenvalue}
	\scalebox{0.58}{
	\begin{tabular}{lccccccc}
	\hline
			$\epsilon$     &       & 1     & $10^{-1}$   & $10^{-2}$  & $10^{-3}$ & $10^{-4}$ & $10^{-5}$ \\
		\hline
    $VFeSb$ & TRUE  & 4.87650627  & 4.30866681  & 4.25706047  & 4.25195125  & 4.25144085  & 4.25138981  \\
          & (2.1) & 1.35542725  & 3.96178845  & 4.22242458  & 4.24848819  & 4.25109455  & 4.25135518  \\
          & (2.4) & 1.33942214  & 3.96018794  & 4.22226452  & 4.24847218  & 4.25109295  & 4.25135502  \\
          & (2.5) & \textbf{1.61331416 } & \textbf{3.97525310 } & \textbf{4.22371597 } & \textbf{4.24861680 } & \textbf{4.25110740 } & \textbf{4.25135647 } \\
          \hline
    $SiO_2$  & TRUE  & 2.98257552  & 0.39023942  & 0.15777358  & 0.13896998  & 0.13762874  & 0.13754487  \\
          & (2.1) & -2.75842067  & -0.15205946  & 0.10857666  & 0.13464027  & 0.13724663  & 0.13750727  \\
          & (2.4) & -2.77442578  & -0.15365997  & 0.10841661  & 0.13462427  & 0.13724503  & 0.13750711  \\
          & (2.5) & \textbf{0.07971937 } & \textbf{0.09417936 } & \textbf{0.13108013 } & \textbf{0.13686338 } & \textbf{0.13746866 } & \textbf{0.13752947 } \\
    \hline
    $Cr_2AgBiO_8$ & TRUE  & 3.74259460  & 2.70447034  & 2.63169114  & 2.62632963  & 2.62585024  & 2.62580287  \\
          & (2.1) & -0.27015928  & 2.33620192  & 2.59683804  & 2.62290165  & 2.62550802  & 2.62576865  \\
          & (2.4) & -0.28616440  & 2.33460141  & 2.59667799  & 2.62288565  & 2.62550642  & 2.62576849  \\
          & (2.5) & \textbf{1.61256461 } & \textbf{2.47855084 } & \textbf{2.61061075 } & \textbf{2.62427433 } & \textbf{2.62564524 } & \textbf{2.62578237 } \\
    \hline
    $RbTaO_3$ & TRUE  & 13.01446669  & 13.56705411  & 13.63089784  & 13.63738109  & 13.63803042  & 13.63809536  \\
          & (2.1) & 10.74214568  & 13.34850689  & 13.60914301  & 13.63520662  & 13.63781298  & 13.63807362  \\
          & (2.4) & 10.72614057  & 13.34690638  & 13.60898296  & 13.63519062  & 13.63781138  & 13.63807346  \\
          & (2.5) & \textbf{12.02344640 } & \textbf{13.46425668 } & \textbf{13.62060290 } & \textbf{13.63635147 } & \textbf{13.63792745 } & \textbf{13.63808506 } \\
    \hline
    $NaBiS_2$ & TRUE  & 11.57301393  & 11.65711348  & 11.66633731  & 11.66726751  & 11.66736061  & 11.66736992  \\
          & (2.1) & 8.77141406  & 11.37777527  & 11.63841139  & 11.66447500  & 11.66708136  & 11.66734200  \\
          & (2.4) & 8.75540895  & 11.37617476  & 11.63825134  & 11.66445899  & 11.66707976  & 11.66734184  \\
          & (2.5) & \textbf{10.49111918 } & \textbf{11.52865892 } & \textbf{11.65330277 } & \textbf{11.66596218 } & \textbf{11.66723006 } & \textbf{11.66735687 } \\
          \hline
    $LiBiB_2O_5$ & TRUE  & 9.48563900  & 7.88963674  & 7.75257003  & 7.73911060  & 7.73776715  & 7.73763283  \\
          & (2.1) & 4.84166101  & 7.44802222  & 7.70865834  & 7.73472195  & 7.73732831  & 7.73758895  \\
          & (2.4) & 4.82565590  & 7.44642171  & 7.70849829  & 7.73470595  & 7.73732671  & 7.73758879  \\
          & (2.5) & \textbf{7.61350030 } & \textbf{7.71767528 } & \textbf{7.73554134 } & \textbf{7.73740942 } & \textbf{7.73759705 } & \textbf{7.73761582 } \\
          \hline
    $KBi_2F_7$ & TRUE  & 14.68094239  & 13.60467138  & 13.51224483  & 13.50314448  & 13.50223585  & 13.50214501  \\
          & (2.1) & 10.60617802  & 13.21253922  & 13.47317534  & 13.49923895  & 13.50184532  & 13.50210595  \\
          & (2.4) & 10.59017290  & 13.21093871  & 13.47301529  & 13.49922295  & 13.50184372  & 13.50210579  \\
          & (2.5) & \textbf{13.41941345 } & \textbf{13.49366946 } & \textbf{13.50128625 } & \textbf{13.50205002 } & \textbf{13.50212642 } & \textbf{13.50213406 } \\
    \hline
    $BaNiO_3$ & TRUE  & 27.54177747  & 27.46666969  & 27.46314869  & 27.46283449  & 27.46280345  & 27.46280034  \\
          & (2.1) & 24.56684311  & 27.17320431  & 27.43384043  & 27.45990404  & 27.46251040  & 27.46277104  \\
          & (2.4) & 24.55083799  & 27.17160380  & 27.43368038  & 27.45988804  & 27.46250880  & 27.46277088  \\
          & (2.5) & \textbf{26.99366189 } & \textbf{27.41524501 } & \textbf{27.45803818 } & \textbf{27.46232376 } & \textbf{27.46275237 } & \textbf{27.46279524 } \\
    \hline
    \end{tabular}}%
  \label{tab_low}%
\end{table}%
\par
From the bounds in Table 1 and Table 2, we see that both the upper bounds and the lower bounds of each example are close to the true $C$-eigenvalues when the norm of the perturbation tensors becomes small, and the bounds given by (\ref{th2.2-01}) perform better than the others.
\clearpage
\textbf{Declarations}
\\
\textbf{Funding}
The first author was funded by the Science and Technology Research Program of Chongqing Municipal Education Commission (Grant Nos. KJQN202100505, KJQN202200512), the Natural Science Foundation Project of Chongqing (Grant Nos. cstc2021jcyj-msxmX0195, CSTB2022NSCQ-MSX0896). The third author was funded by National Natural Science Foundation of China (No. 12101136), Guangdong Basic and Applied Basic Research Foundations (Nos.  2023A1515011633, 2020A1515110967), Project of Science and Technology of Guangzhou (No. 202102020273), the Open Project of Key Laboratory, School of Mathematical Sciences, Chongqing Normal University (No. CSSXKFKTQ202002).
\\
\textbf{Data availability}	Data sharing not applicable to this article as no datasets were generated or analyzed during the current study.
\\
\textbf{Conflict of interest} The authors declare no competing interests.
\\
\textbf{Consent for publication} All authors reviewed the manuscript and agree to publish it.
\\
\textbf{Acknowledgments} All authors express gratitude to the editors and reviewers for handling this manuscript.
\\
\textbf{Ethical Approval and Consent to participate} Not Applicable.
\\
\textbf{Human and Animal Ethics} Not Applicable.


\end{document}